# ESTIMATION OF WEIBULL SHAPE PARAMETER BY SHRINKAGE TOWARDS AN INTERVAL UNDER FAILURE CENSORED SAMPLING

Housila P. Singh[1], Sharad Saxena[1], Jack Allen[2], Sarjinder Singh[3] and Florentin Smarandache [4]

[1] School of Studies in Statistics, Vikram University, Ujjain - 456 010 (M. P.), India.
[2] School of Accounting and Finance, Griffith University, Australia
[3] Department of Mathematics and Statistics, University of Saskatchewan, Canada
[4] Department of Mathematics, University of New Mexico, USA

*Abstract*

This paper is speculated to propose a class of shrinkage estimators for shape parameter $\beta$ in failure censored samples from two-parameter Weibull distribution when some 'apriori' or guessed interval containing the parameter $\beta$ is available in addition to sample information and analyses their properties. Some estimators are generated from the proposed class and compared with the minimum mean squared error (MMSE) estimator. Numerical computations in terms of percent relative efficiency and absolute relative bias indicate that certain of these estimators substantially improve the MMSE estimator in some guessed interval of the parameter space of $\beta$, especially for censored samples with small sizes. Subsequently, a modified class of shrinkage estimators is proposed with its properties.

## 1. INTRODUCTION

Identical rudiments subjected to identical environmental conditions will fail at different and unpredictable times. The 'time of failure' or 'life length' of a component, measured from some specified time until it fails, is represented by the continuous random variable *X*. One distribution that has been used extensively in recent years to deal with such problems of reliability and life-testing is the Weibull distribution introduced by Weibull(1939), who proposed it in connection with his studies on strength of material.

The Weibull distribution includes the exponential and the Rayleigh distributions as special cases. The use of the distribution in reliability and quality control work was advocated by many authors following Weibull(1951), Lieblin and Zelen(1956), Kao(1958,1959), Berrettoni(1964) and Mann(1968 A). Weibull(1951) showed that the distribution is useful in describing the 'wear-out' or fatigue failures. Kao(1959) used it as a model for vacuum tube failures while Lieblin and Zelen(1956) used it as a model for ball bearing failures. Mann(1968

A) gives a variety of situations in which the distribution is used for other types of failure data. The distribution often becomes suitable where the conditions for "strict randomness" of the exponential distribution are not satisfied with the shape parameter $\beta$ having a characteristic or predictable value depending upon the fundamental nature of the problem being considered.

**1.1 The Model**

Let $x_1, x_2, \ldots, x_n$ be a random sample of size $n$ from a two-parameter Weibull distribution, probability density function of which is given by :

$$f(x;\alpha,\beta) = \beta \alpha^{-\beta} x^{\beta-1} \exp\{-(x/\alpha)^\beta\}; \ x > 0, \alpha > 0, \beta > 0 \tag{1.1}$$

where $\alpha$ being the characteristic life acts as a scale parameter and $\beta$ is the shape parameter.

The variable $Y = \ln x$ follows an extreme value distribution, sometimes called the log-Weibull distribution [e.g. White(1969)], cumulative distribution function of which is given by :

$$F(y) = 1 - \exp\left\{-\exp\left(\frac{y-u}{b}\right)\right\}; \ -\infty < y < \infty, -\infty < u < \infty, b > 0 \tag{1.2}$$

where $b = 1/\beta$ and $u = \ln \alpha$ are respectively the scale and location parameters.

The inferential procedures of the above model are quite complex. Mann(1967 A,B, 1968 B) suggested the generalised least squares estimator using the variances and covariances of the ordered observations for which tables are available up to $n = 25$ only.

**1.2 Classical Estimators**

Suppose $x_1, x_2, \ldots, x_m$ be the $m$ smallest ordered observations in a sample of size $n$ from Weibull distribution. Bain(1972) defined an unbiased estimator for $b$ as

$$\hat{b}_u = -\sum_{i=1}^{m-1} \left\{\frac{y_i - y_m}{nK_{(m,n)}}\right\}, \tag{1.3}$$

where 
$$K_{(m,n)} = -\left(\frac{1}{n}\right) E\left[\sum_{i=1}^{m-1}(v_i - v_m)\right], \tag{1.4}$$

and $v_i = \dfrac{y_i - u}{b}$ are ordered variables from the extreme value distribution with $u = 0$ and $b = 1$. The estimator $\hat{b}_u$ is found to have high relative efficiency for heavily censored cases. Contrary to this, the asymptotic relative efficiency of $\hat{b}_u$ is zero for complete samples.

Engelhardt and Bain(1973) suggested a general form of the estimator as

$$\hat{b}_g = -\sum_{i=1}^{m} \left\{ \frac{|y_i - y_m|}{nK_{(g,m,n)}} \right\}, \tag{1.5}$$

where $g$ is a constant to be chosen so that the variance of $\hat{b}_g$ is least and $K_{(g,m,n)}$ is an unbiasing constant. The statistic $\dfrac{h\hat{b}_g}{b}$ has been shown to follow approximately $\chi^2$ - distribution with $h$ degrees of freedom, where $h = 2 \big/ Var(\hat{b}_g/b)$. Therefore, we have

$$E\left[\hat{\beta}^{-jp}\right] = \frac{1}{\beta^{jp}} \left(\frac{2}{h-2}\right)^{jp} \frac{\Gamma[(h/2)+jp]}{\Gamma(h/2)} \; ; \; j = 1,2 \tag{1.6}$$

where $\hat{\beta} = \dfrac{h-2}{t}$ is an unbiased estimator of $\beta$ with $Var(\hat{\beta}) = \dfrac{2\beta^2}{(h-4)}$ and $t = h\hat{b}_g$ having density

$$f(t) = \frac{1}{\Gamma(h/2)} \left(\frac{\beta}{2}\right)^{h/2} \exp\left(\frac{-\beta t}{2}\right) t^{(h/2)-1} \; ; \; t > 0.$$

The MMSE estimator of $\beta$, among the class of estimators of the form $C\hat{\beta}$ ; $C$ being a constant for which the mean square error (MSE) of $C\hat{\beta}$ is minimum, is

$$\hat{\beta}_M = \frac{h-4}{t}, \tag{1.7}$$

having absolute relative bias and relative mean squared error as

$$\text{ARB}\{\hat{\beta}_M\} = \left|\frac{2}{h-2}\right|, \tag{1.8}$$

and

$$\text{RMSE}\{\hat{\beta}_M\} = \frac{2}{h-2}, \tag{1.9}$$

respectively.

## 1.3 Shrinkage Technique of Estimation

Considerable amount of work dealing with shrinkage estimation methods for the parameters of the Weibull distribution has been done since 1970. An experimenter involved in life-testing experiments becomes quite familiar with failure data and hence may often develop knowledge about some parameters of the distribution. In the case of Weibull distribution, for example, knowledge on the shape parameter $\beta$ can be utilised to develop improved inference for the other parameters. Thompson(1968 A,B) considered the problem of shrinking an unbiased estimator $\hat{\xi}$ of the parameter $\xi$ either towards a natural origin $\xi_0$ or towards an interval $(\xi_1, \xi_2)$ and suggested the shrunken estimators $h\hat{\xi} + (1-h)\xi_0$ and $h\hat{\xi} + (1-h)\left(\dfrac{\xi_1 + \xi_2}{2}\right)$, where $0 < h < 1$ is a constant. The relevance of such type of shrunken estimators lies in the fact that, though perhaps they are biased, has smaller MSE than $\hat{\xi}$ for $\xi$ in some interval around $\xi_0$ or $\left(\dfrac{\xi_1 + \xi_2}{2}\right)$, as the case may be. This type of shrinkage estimation of the Weibull parameters has been discussed by various authors, including Singh and Bhatkulikar(1978), Pandey(1983), Pandey and Upadhyay(1985,1986) and Singh and Shukla(2000). For example, Singh and Bhatkulikar(1978) suggested performing a significance test of the validity of the prior value of $\beta$ (which they took as 1). Pandey(1983) also suggested a similar preliminary test shrunken estimator for $\beta$.

In the present investigation, it is desired to estimate $\beta$ in the presence of a prior information available in the form of an interval $(\beta_1, \beta_2)$ and the sample information contained in $\hat{\beta}$. Consequently, this article is an attempt in the direction of obtaining an efficient class of shrunken estimators for the scale parameter $\beta$. The properties of the suggested class of estimators are also discussed theoretically and empirically. The proposed class of shrunken estimators is furthermore modified with its properties.

## 2. THE PROPOSED CLASS OF SHRINKAGE ESTIMATORS

Consider a class of estimators $\overset{*}{\beta}_{(p,q)}$ for $\beta$ in model (1.1) defined by

$$\overset{*}{\beta}_{(p,q)} = \left(\frac{\beta_1 + \beta_2}{2}\right)\left[q + w\left(\frac{\beta_1 + \beta_2}{2\hat{\beta}}\right)^p\right], \qquad (2.1)$$

where $p$ and $q$ are real numbers such that $p \neq 0$ and $q > 0$, $w$ is a stochastic variable which may in particular be a scalar, to be chosen such that MSE of $\overset{*}{\beta}_{(p,q)}$ is minimum.

Assuming $w$ a scalar and using result (1.6), the MSE of $\overset{*}{\beta}_{(p,q)}$ is given by

$$\text{MSE}\left\{\overset{*}{\beta}_{(p,q)}\right\} = \beta^2\left[\{q\Delta-1\}^2 + w^2\Delta^{2(p+1)}\left(\frac{2}{h-2}\right)^{2p}\frac{\Gamma[(h/2)+2p]}{\Gamma(h/2)}\right.$$

$$\left.+\{q\Delta-1\}w\Delta^{(p+1)}\left(\frac{2}{h-2}\right)^{p}\frac{\Gamma[(h/2)+2p]}{\Gamma(h/2)}\right] \quad (2.2)$$

where $\quad \Delta = \left(\dfrac{\beta_1+\beta_2}{2\beta}\right)$.

Minimising (2.2) with respect to $w$ and replacing $\beta$ by its unbiased estimator $\hat{\beta}$, we get

$$\hat{w} = \frac{-\left\{q\left(\dfrac{\beta_1+\beta_2}{2}\right)-\hat{\beta}\right\}\hat{\beta}^p}{\left(\dfrac{\beta_1+\beta_2}{2}\right)^{(p+1)}}w(p). \quad (2.3)$$

where $\quad w(p) = \left(\dfrac{h-2}{2}\right)^p \dfrac{\Gamma[(h/2)+p]}{\Gamma[(h/2)+2p]},$ \quad (2.4)

lies between 0 and 1, {i.e., $0 < w(p) \leq 1$} provided gamma functions exist, i.e., $p > (-h/2)$.

Substituting (2.3) in (2.1) yields a class of shrinkage estimators for $\beta$ in a more feasible form as

$$\hat{\beta}_{(p,q)} = \left(\frac{h-2}{t}\right)w(p) + q\left(\frac{\beta_1+\beta_2}{2}\right)\{1-w(p)\}. \quad (2.5)$$

**2.1 Non-negativity**

Clearly, the proposed class of estimators (2.5) is the convex combination of $\{(h-2)/t\}$ and $\{q(\beta_1+\beta_2)/2\}$ and hence $\hat{\beta}_{(p,q)}$ is always positive as $\{(h-2)/t\} > 0$ and $q > 0$.

**2.2 Unbiasedness**

If $w(p) = 1$, the proposed class of shrinkage estimators $\hat{\beta}_{(p,q)}$ turns into the unbiased estimator $\hat{\beta}$, otherwise it is biased with

$$\text{Bias}\left\{\hat{\beta}_{(p,q)}\right\} = \beta\{q\Delta-1\}[1-w(p)] \quad (2.6)$$

and thus the absolute relative bias is given by

$$\text{ARB}\left\{\hat{\beta}_{(p,q)}\right\} = \left|\{q\Delta - 1\}[1 - w(p)]\right|. \tag{2.7}$$

The condition for unbiasedness that $w(p) = 1$, holds iff, censored sample size $m$ is indefinitely large, i.e., $m \to \infty$. Moreover, if the proposed class of estimators $\hat{\beta}_{(p,q)}$ turns into $\hat{\beta}$ then this case does not deal with the use of prior information.

A more realistic condition for unbiasedness without damaging the basic structure of $\hat{\beta}_{(p,q)}$ and utilises prior information intelligibly can be obtained by (2.7). The ARB of $\hat{\beta}_{(p,q)}$ is zero when $q = \Delta^{-1}$ (or $\Delta = q^{-1}$).

**2.3 Relative Mean Squared Error**

The MSE of the suggested class of shrinkage estimators is derived as

$$\text{MSE}\left\{\hat{\beta}_{(p,q)}\right\} = \beta^2\left[\{q\Delta - 1\}^2\{1 - w(p)\}^2 + \frac{2\{w(p)\}^2}{(h-4)}\right], \tag{2.8}$$

and relative mean square error is therefore given by

$$\text{RMSE}\left\{\hat{\beta}_{(p,q)}\right\} = \{q\Delta - 1\}^2\{1 - w(p)\}^2 + \frac{2\{w(p)\}^2}{(h-4)}. \tag{2.9}$$

It is obvious from (2.9) that $\text{RMSE}\left\{\hat{\beta}_{(p,q)}\right\}$ is minimum when $q = \Delta^{-1}$ (or $\Delta = q^{-1}$).

**2.4 Selection of the Scalar '$p$'**

The convex nature of the proposed statistic and the condition that gamma functions contained in $w(p)$ exist, provides the criterion of choosing the scalar $p$. Therefore, the acceptable range of value of $p$ is given by

$$\{p \mid 0 < w(p) \leq 1 \text{ and } p > (-h/2)\}, \; \forall \; n, m. \tag{2.10}$$

**2.5 Selection of the Scalar '$q$'**

It is pointed out that at $q = \Delta^{-1}$, the proposed class of estimators is not only unbiased but renders maximum gain in efficiency, which is a remarkable property of the proposed class of estimators. Thus obtaining significant gain in efficiency as well as proportionately small magnitude of bias for fixed $\Delta$ or for fixed $(\beta_1/\beta)$ and $(\beta_2/\beta)$, one should choose $q$ in the vicinity of $q = \Delta^{-1}$. It is interesting to note that if one selects smaller values of $q$ then higher values of $\Delta$ leads to a large gain in efficiency (along with appreciable smaller magnitude of bias) and vice-versa. This implies that for smaller values of $q$, the proposed class of estimators allows to choose the guessed interval much wider, i.e., even if the experimenter is less experienced the risk of estimation using the proposed class of estimators is not higher. This is legitimate for all values of $p$.

## 2.3 Estimation of Average Departure: A Practical Way of selecting $q$

The quantity $\Delta = \{(\beta_1 + \beta_2)/2\beta\}$, represents the average departure of natural origins $\beta_1$ and $\beta_2$ from the true value $\beta$. But in practical situations it is hardly possible to get an idea about $\Delta$. Consequently, an unbiased estimator of $\Delta$ is proposed, namely

$$\hat{\Delta} = \left\{\frac{t(\beta_1 + \beta_2)}{4}\right\} \frac{\Gamma(h/2)}{\Gamma[(h/2)+1]}. \tag{2.12}$$

In section 2.5 it is investigated that, if $q = \Delta^{-1}$, the suggested class of estimators yields favourable results. Keeping in view of this concept, one may select $q$ as

$$q = \hat{\Delta}^{-1} = \left\{\frac{4}{t(\beta_1 + \beta_2)}\right\} \frac{\Gamma[(h/2)+1]}{\Gamma(h/2)}. \tag{2.13}$$

Here this is fit for being quoted that this is the criterion of selecting $q$ numerically and one should carefully notice that this doesn't mean $q$ is replaced by (2.13) in $\hat{\beta}_{(p,q)}$.

## 3. COMPARISION OF ESTIMATORS AND EMPIRICAL STUDY

James and Stein(1961) reported that minimum MSE is a highly desirable property and it is therefore used as a criterion to compare different estimators with each other. The condition under which the proposed class of estimators is more efficient than the MMSE estimator is given below.

$\text{MSE}\left\{\hat{\beta}_{(p,q)}\right\}$ does not exceed the MSE of MMSE estimator $\hat{\beta}_M$ if -

$$\left(1 - \sqrt{G}\right)q^{-1} < \Delta < \left(1 + \sqrt{G}\right)q^{-1} \tag{3.1}$$

where $\quad G = \dfrac{2}{\{1 - w(p)\}^2}\left[\dfrac{1}{(h-2)} - \dfrac{\{w(p)\}^2}{(h-4)}\right].$

Besides minimum MSE criterion, minimum bias is also important and therefore should be incorporated under study. Thus, ARB$\{\hat{\beta}_{(p,q)}\}$ is less than ARB$\{\hat{\beta}_M\}$ if -

$$\left\{1-\frac{2}{(h-2)(1-w_{(p)})}\right\}q^{-1} < \Delta < \left\{1+\frac{2}{(h-2)(1-w_{(p)})}\right\}q^{-1} \quad (3.2)$$

### 3.1 The Best Range of Dominance of $\Delta$

The intersection of the ranges of $\Delta$ in (3.1) and (3.2) gives the best range of dominance of $\Delta$ denoted by $\Delta_{Best}$. In this range, the proposed class of estimators is not only less biased than the MMSE estimator but is more efficient than that. The four possible cases in this regard are:

(i) if $\left\{1-\frac{2}{(h-2)[1-w(p)]}\right\} < (1-\sqrt{G})$ and $\left\{1+\frac{2}{(h-2)[1-w(p)]}\right\} < (1+\sqrt{G})$ then

$$\Delta_{Best} = \left(\{1-\sqrt{G}\}q^{-1}, \left\{1+\frac{2}{(h-2)[1-w(p)]}\right\}q^{-1}\right)$$

(ii) if $\left\{1-\frac{2}{(h-2)[1-w(p)]}\right\} < (1-\sqrt{G})$ and $(1+\sqrt{G}) < \left\{1+\frac{2}{(h-2)[1-w(p)]}\right\}$ then

$\Delta_{Best}$ is the same as defined in (3.1).

(iii) if $(1-\sqrt{G}) < \left\{1-\frac{2}{(h-2)[1-w(p)]}\right\}$ and $(1+\sqrt{G}) < \left\{1+\frac{2}{(h-2)[1-w(p)]}\right\}$ then

$$\Delta_{Best} = \left(\left\{1-\frac{2}{(h-2)[1-w(p)]}\right\}q^{-1}, \{1+\sqrt{G}\}q^{-1}\right)$$

(iv) if $(1-\sqrt{G}) < \left\{1-\frac{2}{(h-2)[1-w(p)]}\right\}$ and $\left\{1+\frac{2}{(h-2)[1-w(p)]}\right\} < (1+\sqrt{G})$ then

$\Delta_{Best}$ is the same as defined in (3.2).

### 3.2 Percent Relative Efficiency

To elucidate the performance of the proposed class of estimators $\hat{\beta}_{(p,q)}$ with the MMSE estimator $\hat{\beta}_M$, the Percent Relative Efficiencies (PREs) of $\hat{\beta}_{(p,q)}$ with respect to $\hat{\beta}_M$ have been computed by the formula:

$$\text{PRE}\left\{\hat{\beta}_{(p,q)}, \hat{\beta}_M\right\} = \frac{2(h-4)}{(h-2)\left[(q\Delta-1)^2\{1-w(p)\}^2(h-4)+2\{w(p)\}^2\right]} \times 100 \qquad (3.5)$$

The PREs of $\hat{\beta}_{(p,q)}$ with respect to $\hat{\beta}_M$ and ARBs of $\hat{\beta}_{(p,q)}$ for fixed $n = 20$ and different values of $p$, $q$, $m$ $\Delta_1(=\beta_1/\beta)$ and $\Delta_2(=\beta_2/\beta)$ or $\Delta$ are compiled in Table 3.1 with corresponding values of $h$ [which can be had from Engelhardt(1975)] and $w(p)$. The first column in every $m$ corresponds to PREs and the second one corresponds to ARBs of $\hat{\beta}_{(p,q)}$. The last two rows of each set of $q$ includes the range of dominance of $\Delta$ and $\Delta_{Best}$. The ARBs of $\hat{\beta}_M$ has also been given at the end of each set of table.

**Table 3.1**

PREs of proposed estimator $\hat{\hat{\beta}}_{(p,q)}$ with respect to MMSE estimator $\hat{\hat{\beta}}_m$ and ARBs of $\hat{\hat{\beta}}_{(p,q)}$

| $q\downarrow$ | $\Delta_1\downarrow$ | $\Delta_2\downarrow$ | $m\rightarrow$ $h\rightarrow$ $\Delta\downarrow$ $w(p)\rightarrow$ | \multicolumn{2}{c}{$p = -2$} | | | | | | |
|---|---|---|---|---|---|---|---|---|---|---|---|
| | | | | \multicolumn{2}{c}{6} | \multicolumn{2}{c}{8} | \multicolumn{2}{c}{10} | \multicolumn{2}{c}{12} |
| | | | | \multicolumn{2}{c}{10.8519} | \multicolumn{2}{c}{15.6740} | \multicolumn{2}{c}{20.8442} | \multicolumn{2}{c}{26.4026} |
| | | | | \multicolumn{2}{c}{0.1750} | \multicolumn{2}{c}{0.3970} | \multicolumn{2}{c}{0.5369} | \multicolumn{2}{c}{0.6305} |
| 0.25 | 0.1 | 0.2 | 0.15 | 35.33 | 0.7941 | 40.20 | 0.5804 | 45.57 | 0.4457 | 50.60 | 0.3556 |
| | 0.4 | 0.6 | 0.50 | 42.62 | 0.7219 | 47.90 | 0.5276 | 53.49 | 0.4052 | 58.53 | 0.3233 |
| | 0.4 | 1.6 | 1.00 | 57.66 | 0.6188 | 63.18 | 0.4522 | 68.54 | 0.3473 | 72.99 | 0.2771 |
| | 1.0 | 2.0 | 1.50 | 82.21 | 0.5156 | 86.53 | 0.3769 | 89.95 | 0.2894 | 92.27 | 0.2309 |
| | 1.6 | 2.4 | 2.00 | 126.15 | 0.4125 | 124.06 | 0.3015 | 120.83 | 0.2315 | 117.72 | 0.1847 |
| | 2.0 | 3.0 | 2.50 | 215.89 | 0.3094 | 187.20 | 0.2261 | 164.84 | 0.1737 | 149.86 | 0.1386 |
| | 2.5 | 3.5 | 3.00 | 438.90 | 0.2063 | 294.12 | 0.1507 | 222.82 | 0.1158 | 186.17 | 0.0924 |
| | 3.5 | 3.5 | 3.50 | 1154.45 | 0.1031 | 447.47 | 0.0754 | 282.42 | 0.0579 | 217.84 | 0.0462 |
| | 3.8 | 4.2 | 4.00 | 2528.52 | 0.0000 | 541.60 | 0.0000 | 310.07 | 0.0000 | 230.93 | 0.0000 |
| | \multicolumn{3}{c}{Range of $\Delta\rightarrow$} | (1.74, 6.25) | (2.90, 5.09) | (1.70, 6.29) | (3.02, 4.97) | (1.68, 6.31) | (3.08, 4.91) | (1.66, 6.33) | (3.11, 4.88) |
| | \multicolumn{3}{c}{$\Delta_{Best}\rightarrow$} | \multicolumn{2}{c}{(2.90, 5.09)} | \multicolumn{2}{c}{(3.02, 4.97)} | \multicolumn{2}{c}{(3.08, 4.91)} | \multicolumn{2}{c}{(3.11, 4.88)} |
| 0.50 | 0.1 | 0.2 | 0.15 | 38.21 | 0.7632 | 43.26 | 0.5577 | 48.75 | 0.4284 | 53.81 | 0.3418 |
| | 0.4 | 0.6 | 0.50 | 57.66 | 0.6188 | 63.18 | 0.4522 | 68.54 | 0.3473 | 72.99 | 0.2771 |
| | 0.4 | 1.6 | 1.00 | 126.15 | 0.4125 | 124.06 | 0.3015 | 120.83 | 0.2315 | 117.72 | 0.1847 |
| | 1.0 | 2.0 | 1.50 | 438.90 | 0.2063 | 294.12 | 0.1507 | 222.82 | 0.1158 | 186.17 | 0.0924 |
| | 1.6 | 2.4 | 2.00 | 2528.52 | 0.0000 | 541.60 | 0.0000 | 310.07 | 0.0000 | 230.93 | 0.0000 |
| | 2.0 | 3.0 | 2.50 | 438.90 | 0.2063 | 294.12 | 0.1507 | 222.82 | 0.1158 | 186.17 | 0.0924 |
| | 2.5 | 3.5 | 3.00 | 126.15 | 0.4125 | 124.06 | 0.3015 | 120.83 | 0.2315 | 117.72 | 0.1847 |
| | 3.5 | 3.5 | 3.50 | 57.66 | 0.6188 | 63.18 | 0.4522 | 68.54 | 0.3473 | 72.99 | 0.2771 |
| | 3.8 | 4.2 | 4.00 | 32.76 | 0.8250 | 37.45 | 0.6030 | 42.68 | 0.4631 | 47.65 | 0.3695 |
| | \multicolumn{3}{c}{Range of $\Delta\rightarrow$} | (0.87, 3.13) | (1.45, 2.55) | (0.85, 3.15) | (1.51, 2.49) | (0.84, 3.16) | (1.54, 2.46) | (0.83, 3.17) | (1.56, 2.44) |
| | \multicolumn{3}{c}{$\Delta_{Best}\rightarrow$} | \multicolumn{2}{c}{(1.45, 2.55)} | \multicolumn{2}{c}{(1.51, 2.49)} | \multicolumn{2}{c}{(1.54, 2.46)} | \multicolumn{2}{c}{(1.56, 2.44)} |
| 0.75 | 0.1 | 0.2 | 0.15 | 41.45 | 0.7322 | 46.67 | 0.5351 | 52.25 | 0.4110 | 57.30 | 0.3279 |
| | 0.4 | 0.6 | 0.50 | 82.21 | 0.5156 | 86.53 | 0.3769 | 89.95 | 0.2894 | 92.27 | 0.2309 |
| | 0.4 | 1.6 | 1.00 | 438.90 | 0.2063 | 294.12 | 0.1507 | 222.82 | 0.1158 | 186.17 | 0.0924 |
| | 1.0 | 2.0 | 1.50 | 1154.45 | 0.1031 | 447.47 | 0.0754 | 282.42 | 0.0579 | 217.84 | 0.0462 |
| | 1.6 | 2.4 | 2.00 | 126.15 | 0.4125 | 124.06 | 0.3015 | 120.83 | 0.2315 | 117.72 | 0.1847 |
| | 2.0 | 3.0 | 2.50 | 42.62 | 0.7219 | 47.90 | 0.5276 | 53.49 | 0.4052 | 58.53 | 0.3233 |
| | 2.5 | 3.5 | 3.00 | 21.07 | 1.0313 | 24.58 | 0.7537 | 28.74 | 0.5789 | 32.94 | 0.4619 |
| | 3.5 | 3.5 | 3.50 | 12.51 | 1.3407 | 14.82 | 0.9798 | 17.67 | 0.7525 | 20.70 | 0.6004 |
| | 3.8 | 4.2 | 4.00 | 8.27 | 1.6501 | 9.87 | 1.2059 | 11.90 | 0.9262 | 14.09 | 0.7390 |
| | \multicolumn{3}{c}{Range of $\Delta\rightarrow$} | (0.58, 2.09) | (0.97, 1.70) | (0.57, 2.10) | (1.01, 1.66) | (0.56, 2.11) | (1.03, 1.64) | (0.56, 2.11) | (1.04, 1.63) |
| | \multicolumn{3}{c}{$\Delta_{Best}\rightarrow$} | \multicolumn{2}{c}{(0.97, 1.70)} | \multicolumn{2}{c}{(1.01, 1.66)} | \multicolumn{2}{c}{(1.03, 1.64)} | \multicolumn{2}{c}{(1.04, 1.63)} |
| \multicolumn{4}{c}{**ARB of MMSE Estimator**$\rightarrow$} | \multicolumn{2}{c}{0.2259} | \multicolumn{2}{c}{0.1463} | \multicolumn{2}{c}{0.1061} | \multicolumn{2}{c}{0.0820} |

*Table 3.1 continued …*

| $q\downarrow$ | $\Delta_1\downarrow$ | $\Delta_2\downarrow$ | $m\rightarrow$ | \multicolumn{8}{c}{$p = -1$} |||||||||
|---|---|---|---|---|---|---|---|---|---|---|---|
| | | | | \multicolumn{2}{c\|}{6} | \multicolumn{2}{c\|}{8} | \multicolumn{2}{c\|}{10} | \multicolumn{2}{c}{12} |
| | | | $h\rightarrow$ | \multicolumn{2}{c\|}{10.8519} | \multicolumn{2}{c\|}{15.6740} | \multicolumn{2}{c\|}{20.8442} | \multicolumn{2}{c}{26.4026} |
| | | | $\Delta\downarrow\ w(p)\rightarrow$ | \multicolumn{2}{c\|}{0.7739} | \multicolumn{2}{c\|}{0.8537} | \multicolumn{2}{c\|}{0.8939} | \multicolumn{2}{c}{0.9180} |
| 0.25 | 0.1 | 0.2 | 0.15 | 101.69 | 0.2176 | 101.09 | 0.1408 | 100.79 | 0.1022 | 100.61 | 0.0789 |
| | 0.4 | 0.6 | 0.50 | 105.60 | 0.1978 | 103.55 | 0.1280 | 102.55 | 0.0929 | 101.96 | 0.0718 |
| | 0.4 | 1.6 | 1.00 | 110.98 | 0.1696 | 106.84 | 0.1097 | 104.87 | 0.0796 | 103.73 | 0.0615 |
| | 1.0 | 2.0 | 1.50 | 115.99 | 0.1413 | 109.79 | 0.0914 | 106.91 | 0.0663 | 105.27 | 0.0513 |
| | 1.6 | 2.4 | 2.00 | 120.43 | 0.1130 | 112.32 | 0.0731 | 108.65 | 0.0531 | 106.56 | 0.0410 |
| | 2.0 | 3.0 | 2.50 | 124.13 | 0.0848 | 114.38 | 0.0549 | 110.04 | 0.0398 | 107.59 | 0.0308 |
| | 2.5 | 3.5 | 3.00 | 126.91 | 0.0565 | 115.89 | 0.0366 | 111.05 | 0.0265 | 108.34 | 0.0205 |
| | 3.5 | 3.5 | 3.50 | 128.65 | 0.0283 | 116.82 | 0.0183 | 111.67 | 0.0133 | 108.79 | 0.0103 |
| | 3.8 | 4.2 | 4.00 | 129.23 | 0.0000 | 117.13 | 0.0000 | 111.87 | 0.0000 | 108.94 | 0.0000 |
| | \multicolumn{3}{c\|}{**Range of $\Delta\rightarrow$**} | (0.00, 8.00) | (0.00, 8.00) | (0.00, 8.00) | (0.00, 8.00) | (0.00, 8.00) | (0.00, 8.00) | (0.00, 8.00) | (0.00, 8.00) |
| | \multicolumn{3}{c\|}{$\Delta_{Best}\rightarrow$} | \multicolumn{2}{c\|}{(0.00, 8.00)} | \multicolumn{2}{c\|}{(0.00, 8.00)} | \multicolumn{2}{c\|}{(0.00, 8.00)} | \multicolumn{2}{c}{(0.00, 8.00)} |
| 0.50 | 0.1 | 0.2 | 0.15 | 103.38 | 0.2091 | 102.16 | 0.1353 | 101.56 | 0.0982 | 101.20 | 0.0759 |
| | 0.4 | 0.6 | 0.50 | 110.98 | 0.1696 | 106.84 | 0.1097 | 104.87 | 0.0796 | 103.73 | 0.0615 |
| | 0.4 | 1.6 | 1.00 | 120.43 | 0.1130 | 112.32 | 0.0731 | 108.65 | 0.0531 | 106.56 | 0.0410 |
| | 1.0 | 2.0 | 1.50 | 126.91 | 0.0565 | 115.89 | 0.0366 | 111.05 | 0.0265 | 108.34 | 0.0205 |
| | 1.6 | 2.4 | 2.00 | 129.23 | 0.0000 | 117.13 | 0.0000 | 111.87 | 0.0000 | 108.94 | 0.0000 |
| | 2.0 | 3.0 | 2.50 | 126.91 | 0.0565 | 115.89 | 0.0366 | 111.05 | 0.0265 | 108.34 | 0.0205 |
| | 2.5 | 3.5 | 3.00 | 120.43 | 0.1130 | 112.32 | 0.0731 | 108.65 | 0.0531 | 106.56 | 0.0410 |
| | 3.5 | 3.5 | 3.50 | 110.98 | 0.1696 | 106.84 | 0.1097 | 104.87 | 0.0796 | 103.73 | 0.0615 |
| | 3.8 | 4.2 | 4.00 | 100.00 | 0.2261 | 100.00 | 0.1463 | 100.00 | 0.1061 | 100.00 | 0.0820 |
| | \multicolumn{3}{c\|}{**Range of $\Delta\rightarrow$**} | (0.00, 4.00) | (0.00, 4.00) | (0.00, 4.00) | (0.00, 4.00) | (0.00, 4.00) | (0.00, 4.00) | (0.00, 4.00) | (0.00, 4.00) |
| | \multicolumn{3}{c\|}{$\Delta_{Best}\rightarrow$} | \multicolumn{2}{c\|}{(0.00, 4.00)} | \multicolumn{2}{c\|}{(0.00, 4.00)} | \multicolumn{2}{c\|}{(0.00, 4.00)} | \multicolumn{2}{c}{(0.00, 4.00)} |
| 0.75 | 0.1 | 0.2 | 0.15 | 105.05 | 0.2006 | 103.21 | 0.1298 | 102.31 | 0.0942 | 101.77 | 0.0728 |
| | 0.4 | 0.6 | 0.50 | 115.99 | 0.1413 | 109.79 | 0.0914 | 106.91 | 0.0663 | 105.27 | 0.0513 |
| | 0.4 | 1.6 | 1.00 | 126.91 | 0.0565 | 115.89 | 0.0366 | 111.05 | 0.0265 | 108.34 | 0.0205 |
| | 1.0 | 2.0 | 1.50 | 128.65 | 0.0283 | 116.82 | 0.0183 | 111.67 | 0.0133 | 108.79 | 0.0103 |
| | 1.6 | 2.4 | 2.00 | 120.43 | 0.1130 | 112.32 | 0.0731 | 108.65 | 0.0531 | 106.56 | 0.0410 |
| | 2.0 | 3.0 | 2.50 | 105.60 | 0.1978 | 103.55 | 0.1280 | 102.55 | 0.0929 | 101.96 | 0.0718 |
| | 2.5 | 3.5 | 3.00 | 88.71 | 0.2826 | 92.40 | 0.1828 | 94.37 | 0.1327 | 95.59 | 0.1025 |
| | 3.5 | 3.5 | 3.50 | 72.93 | 0.3674 | 80.65 | 0.2377 | 85.17 | 0.1725 | 88.13 | 0.1333 |
| | 3.8 | 4.2 | 4.00 | 59.57 | 0.4521 | 69.50 | 0.2925 | 75.85 | 0.2123 | 80.24 | 0.1640 |
| | \multicolumn{3}{c\|}{**Range of $\Delta\rightarrow$**} | (0.00, 2.67) | (0.00, 2.67) | (0.00, 2.67) | (0.00, 2.67) | (0.00, 2.67) | (0.00, 2.67) | (0.00, 2.67) | (0.00, 2.67) |
| | \multicolumn{3}{c\|}{$\Delta_{Best}\rightarrow$} | \multicolumn{2}{c\|}{(0.00, 2.67)} | \multicolumn{2}{c\|}{(0.00, 2.67)} | \multicolumn{2}{c\|}{(0.00, 2.67)} | \multicolumn{2}{c}{(0.00, 2.67)} |
| \multicolumn{4}{l\|}{**ARB of MMSE Estimator$\rightarrow$**} | \multicolumn{2}{c\|}{0.2259} | \multicolumn{2}{c\|}{0.1463} | \multicolumn{2}{c\|}{0.1061} | \multicolumn{2}{c}{0.0820} |

*Table 3.1 continued …*

| $q\downarrow$ | $\Delta_1\downarrow$ | $\Delta_2\downarrow$ | $m\rightarrow$ | 6 | | 8 | | 10 | | 12 | |
|---|---|---|---|---|---|---|---|---|---|---|---|
| | | | $h\rightarrow$ | 10.8519 | | 15.6740 | | 20.8442 | | 26.4026 | |
| | | | $\Delta\downarrow\ w(p)\rightarrow$ | 0.6888 | | 0.7737 | | 0.8251 | | 0.8779 | |
| | | | | | | | | | | | |
| 0.25 | 0.1 | 0.2 | 0.15 | 99.00 | 0.2996 | 97.51 | 0.2178 | 97.21 | 0.1684 | 99.20 | 0.1175 |
| | 0.4 | 0.6 | 0.50 | 106.26 | 0.2723 | 103.17 | 0.1980 | 101.80 | 0.1531 | 102.17 | 0.1069 |
| | 0.4 | 1.6 | 1.00 | 117.09 | 0.2334 | 111.34 | 0.1697 | 108.25 | 0.1312 | 106.18 | 0.0916 |
| | 1.0 | 2.0 | 1.50 | 128.15 | 0.1945 | 119.34 | 0.1415 | 114.39 | 0.1093 | 109.82 | 0.0763 |
| | 1.6 | 2.4 | 2.00 | 138.88 | 0.1556 | 126.79 | 0.1132 | 119.95 | 0.0875 | 113.00 | 0.0611 |
| | 2.0 | 3.0 | 2.50 | 148.56 | 0.1167 | 133.27 | 0.0849 | 124.67 | 0.0656 | 115.60 | 0.0458 |
| | 2.5 | 3.5 | 3.00 | 156.33 | 0.0778 | 138.31 | 0.0566 | 128.27 | 0.0437 | 117.53 | 0.0305 |
| | 3.5 | 3.5 | 3.50 | 161.41 | 0.0389 | 141.52 | 0.0283 | 130.54 | 0.0219 | 118.72 | 0.0153 |
| | 3.8 | 4.2 | 4.00 | 163.17 | 0.0000 | 142.63 | 0.0000 | 131.31 | 0.0000 | 119.12 | 0.0000 |
| | | | Range of $\Delta\rightarrow$ | (0.20, 7.80) | (0.00, 8.00) | (0.30, 7.70) | (0.00, 8.00) | (0.36, 7.64) | (0.00, 8.00) | (0.24, 7.76) | (0.00, 8.00) |
| | | | | (0.20, 7.80) | | (0.30, 7.70) | | (0.36, 7.64) | | (0.24, 7.76) | |
| 0.50 | 0.1 | 0.2 | 0.15 | 102.07 | 0.2879 | 99.92 | 0.2093 | 99.18 | 0.1618 | 100.49 | 0.1130 |
| | 0.4 | 0.6 | 0.50 | 117.09 | 0.2334 | 111.34 | 0.1697 | 108.25 | 0.1312 | 106.18 | 0.0916 |
| | 0.4 | 1.6 | 1.00 | 138.88 | 0.1556 | 126.79 | 0.1132 | 119.95 | 0.0875 | 113.00 | 0.0611 |
| | 1.0 | 2.0 | 1.50 | 156.33 | 0.0778 | 138.31 | 0.0566 | 128.27 | 0.0437 | 117.53 | 0.0305 |
| | 1.6 | 2.4 | 2.00 | 163.17 | 0.0000 | 142.63 | 0.0000 | 131.31 | 0.0000 | 119.12 | 0.0000 |
| | 2.0 | 3.0 | 2.50 | 156.33 | 0.0778 | 138.31 | 0.0566 | 128.27 | 0.0437 | 117.53 | 0.0305 |
| | 2.5 | 3.5 | 3.00 | 138.88 | 0.1556 | 126.79 | 0.1132 | 119.95 | 0.0875 | 113.00 | 0.0611 |
| | 3.5 | 3.5 | 3.50 | 117.09 | 0.2334 | 111.34 | 0.1697 | 108.25 | 0.1312 | 106.18 | 0.0916 |
| | 3.8 | 4.2 | 4.00 | 96.01 | 0.3112 | 95.12 | 0.2263 | 95.25 | 0.1749 | 97.90 | 0.1221 |
| | | | Range of $\Delta\rightarrow$ | (0.10, 3.90) | (0.55, 3.45) | (0.15, 3.85) | (0.71, 3.29) | (0.18, 3.82) | (0.79, 3.21) | (0.12, 3.88) | (0.66, 3.34) |
| | | | $\Delta_{Best}\rightarrow$ | (0.55, 3.45) | | (0.71, 3.29) | | (0.79, 3.21) | | (0.66, 3.34) | |
| 0.75 | 0.1 | 0.2 | 0.15 | 105.20 | 0.2762 | 102.36 | 0.2009 | 101.15 | 0.1553 | 101.75 | 0.1084 |
| | 0.4 | 0.6 | 0.50 | 128.15 | 0.1945 | 119.34 | 0.1415 | 114.39 | 0.1093 | 109.82 | 0.0763 |
| | 0.4 | 1.6 | 1.00 | 156.33 | 0.0778 | 138.31 | 0.0566 | 128.27 | 0.0437 | 117.53 | 0.0305 |
| | 1.0 | 2.0 | 1.50 | 161.41 | 0.0389 | 141.52 | 0.0283 | 130.54 | 0.0219 | 118.72 | 0.0153 |
| | 1.6 | 2.4 | 2.00 | 138.88 | 0.1556 | 126.79 | 0.1132 | 119.95 | 0.0875 | 113.00 | 0.0611 |
| | 2.0 | 3.0 | 2.50 | 106.26 | 0.2723 | 103.17 | 0.1980 | 101.80 | 0.1531 | 102.17 | 0.1069 |
| | 2.5 | 3.5 | 3.00 | 77.96 | 0.3891 | 80.11 | 0.2829 | 82.50 | 0.2187 | 88.98 | 0.1526 |
| | 3.5 | 3.5 | 3.50 | 57.31 | 0.5058 | 61.51 | 0.3678 | 65.66 | 0.2843 | 75.76 | 0.1984 |
| | 3.8 | 4.2 | 4.00 | 42.96 | 0.6225 | 47.58 | 0.4526 | 52.22 | 0.3499 | 63.80 | 0.2442 |
| | | | Range of $\Delta\rightarrow$ | (0.07, 2.60) | (0.37, 2.30) | (0.10, 2.57) | (0.47, 2.20) | (0.12, 2.55) | (0.52, 2.14) | (0.08, 2.59) | (0.44, 2.23) |
| | | | $\Delta_{Best}\rightarrow$ | (0.37, 2.30) | | (0.47, 2.20) | | (0.52, 2.14) | | (0.44, 2.23) | |
| **ARB of MMSE Estimator**$\rightarrow$ | | | | 0.2259 | | 0.1463 | | 0.1061 | | 0.0820 | |

Table header spans: $p = 1$

*Table 3.1 continued …*

| $q\downarrow$ | $\Delta_1\downarrow$ | $\Delta_2\downarrow$ | $m\rightarrow$ $h\rightarrow$ $\Delta\downarrow$ $w(p)\rightarrow$ | \multicolumn{8}{c}{$p = 2$} | | | | | | | |
|---|---|---|---|---|---|---|---|---|---|---|---|

| $q\downarrow$ | $\Delta_1\downarrow$ | $\Delta_2\downarrow$ | $\Delta\downarrow$ | 6 | | 8 | | 10 | | 12 | |
|---|---|---|---|---|---|---|---|---|---|---|---|
| | | | $h\rightarrow$ | \multicolumn{2}{c}{10.8519} | \multicolumn{2}{c}{15.6740} | \multicolumn{2}{c}{20.8442} | \multicolumn{2}{c}{26.4026} |
| | | | $w(p)\rightarrow$ | \multicolumn{2}{c}{0.3131} | \multicolumn{2}{c}{0.4385} | \multicolumn{2}{c}{0.5392} | \multicolumn{2}{c}{0.6816} |
| 0.25 | 0.1 | 0.2 | 0.15 | 48.51 | 0.6612 | 45.00 | 0.5405 | 45.90 | 0.4435 | 60.53 | 0.3065 |
| | 0.4 | 0.6 | 0.50 | 57.95 | 0.6011 | 53.31 | 0.4913 | 53.85 | 0.4032 | 68.81 | 0.2786 |
| | 0.4 | 1.6 | 1.00 | 76.84 | 0.5152 | 69.55 | 0.4211 | 68.94 | 0.3456 | 83.20 | 0.2388 |
| | 1.0 | 2.0 | 1.50 | 106.11 | 0.4293 | 93.70 | 0.3509 | 90.35 | 0.2880 | 101.08 | 0.1990 |
| | 1.6 | 2.4 | 2.00 | 154.14 | 0.3435 | 130.87 | 0.2808 | 121.15 | 0.2304 | 122.65 | 0.1592 |
| | 2.0 | 3.0 | 2.50 | 237.92 | 0.2576 | 189.27 | 0.2106 | 164.85 | 0.1728 | 147.06 | 0.1194 |
| | 2.5 | 3.5 | 3.00 | 388.87 | 0.1717 | 277.82 | 0.1404 | 222.08 | 0.1152 | 171.43 | 0.0796 |
| | 3.5 | 3.5 | 3.50 | 627.92 | 0.0859 | 386.26 | 0.0702 | 280.49 | 0.0576 | 190.36 | 0.0398 |
| | 3.8 | 4.2 | 4.00 | 789.74 | 0.0000 | 444.03 | 0.0000 | 307.45 | 0.0000 | 197.63 | 0.0000 |
| | \multicolumn{3}{c|}{Range of $\Delta\rightarrow$} | (1.41, 6.59) | (2.68, 5.32) | (1.60, 6.40) | (2.96, 5.04) | (1.68, 6.32) | (3.08, 4.92) | (1.47, 6.53) | (2.97, 5.03) |
| | \multicolumn{3}{c|}{$\Delta_{Best}\rightarrow$} | \multicolumn{2}{c}{(2.68, 5.32)} | \multicolumn{2}{c}{(2.96, 5.04)} | \multicolumn{2}{c}{(3.08, 4.92)} | \multicolumn{2}{c}{(2.97, 5.03)} |
| 0.50 | 0.1 | 0.2 | 0.15 | 52.26 | 0.6354 | 48.32 | 0.5194 | 49.09 | 0.4262 | 63.91 | 0.2946 |
| | 0.4 | 0.6 | 0.50 | 76.84 | 0.5152 | 69.55 | 0.4211 | 68.94 | 0.3456 | 83.20 | 0.2388 |
| | 0.4 | 1.6 | 1.00 | 154.14 | 0.3435 | 130.87 | 0.2808 | 121.15 | 0.2304 | 122.65 | 0.1592 |
| | 1.0 | 2.0 | 1.50 | 388.87 | 0.1717 | 277.82 | 0.1404 | 222.08 | 0.1152 | 171.43 | 0.0796 |
| | 1.6 | 2.4 | 2.00 | 789.74 | 0.0000 | 444.03 | 0.0000 | 307.45 | 0.0000 | 197.63 | 0.0000 |
| | 2.0 | 3.0 | 2.50 | 388.87 | 0.1717 | 277.82 | 0.1404 | 222.08 | 0.1152 | 171.43 | 0.0796 |
| | 2.5 | 3.5 | 3.00 | 154.14 | 0.3435 | 130.87 | 0.2808 | 121.15 | 0.2304 | 122.65 | 0.1592 |
| | 3.5 | 3.5 | 3.50 | 76.84 | 0.5152 | 69.55 | 0.4211 | 68.94 | 0.3456 | 83.20 | 0.2388 |
| | 3.8 | 4.2 | 4.00 | 45.14 | 0.6869 | 42.00 | 0.5615 | 42.99 | 0.4608 | 57.36 | 0.3184 |
| | \multicolumn{3}{c|}{Range of $\Delta\rightarrow$} | (0.71, 3.29) | (1.34, 2.66) | (0.80, 3.20) | (1.48, 2.52) | (0.84, 3.16) | (1.54, 2.46) | (0.74, 3.26) | (1.49, 2.51) |
| | \multicolumn{3}{c|}{$\Delta_{Best}\rightarrow$} | \multicolumn{2}{c}{(1.34, 2.66)} | \multicolumn{2}{c}{(1.48, 2.52)} | \multicolumn{2}{c}{(1.54, 2.46)} | \multicolumn{2}{c}{(1.49, 2.51)} |
| 0.75 | 0.1 | 0.2 | 0.15 | 56.45 | 0.6096 | 52.00 | 0.4983 | 52.60 | 0.4090 | 67.54 | 0.2826 |
| | 0.4 | 0.6 | 0.50 | 106.11 | 0.4293 | 93.70 | 0.3509 | 90.35 | 0.2880 | 101.08 | 0.1990 |
| | 0.4 | 1.6 | 1.00 | 388.87 | 0.1717 | 277.82 | 0.1404 | 222.08 | 0.1152 | 171.43 | 0.0796 |
| | 1.0 | 2.0 | 1.50 | 627.92 | 0.0859 | 386.26 | 0.0702 | 280.49 | 0.0576 | 190.36 | 0.0398 |
| | 1.6 | 2.4 | 2.00 | 154.14 | 0.3435 | 130.87 | 0.2808 | 121.15 | 0.2304 | 122.65 | 0.1592 |
| | 2.0 | 3.0 | 2.50 | 57.95 | 0.6011 | 53.31 | 0.4913 | 53.85 | 0.4032 | 68.81 | 0.2786 |
| | 2.5 | 3.5 | 3.00 | 29.50 | 0.8587 | 27.83 | 0.7019 | 28.97 | 0.5760 | 41.00 | 0.3980 |
| | 3.5 | 3.5 | 3.50 | 17.73 | 1.1163 | 16.90 | 0.9125 | 17.83 | 0.7488 | 26.50 | 0.5175 |
| | 3.8 | 4.2 | 4.00 | 11.79 | 1.3739 | 11.30 | 1.1230 | 12.01 | 0.9216 | 18.33 | 0.6369 |
| | \multicolumn{3}{c|}{Range of $\Delta\rightarrow$} | (0.47, 2.20) | (0.89, 1.77) | (0.53, 2.13) | (0.99, 1.68) | (0.56, 2.11) | (1.03, 1.64) | (0.49, 2.18) | (0.99, 1.68) |
| | \multicolumn{3}{c|}{$\Delta_{Best}\rightarrow$} | \multicolumn{2}{c}{(0.89, 1.77)} | \multicolumn{2}{c}{(0.99, 1.68)} | \multicolumn{2}{c}{(1.03, 1.64)} | \multicolumn{2}{c}{(0.99, 1.68)} |
| \multicolumn{4}{c|}{**ARB of MMSE Estimator**$\rightarrow$} | \multicolumn{2}{c}{0.2259} | \multicolumn{2}{c}{0.1463} | \multicolumn{2}{c}{0.1061} | \multicolumn{2}{c}{0.0820} |

It has been observed from Table 3.1, that on keeping $m$, $p$, $q$ fixed, the relative efficiencies of the proposed class of shrinkage estimators increases up to $\Delta = q^{-1}$, attains its maximum at this point and then decreases symmetrically in magnitude, as $\Delta$ increases in its range of dominance for all $n$, $p$ and $q$. On the other hand, the ARBs of the proposed class of estimators decreases up to $\Delta = q^{-1}$, the estimator becomes unbiased at this point and then ARBs increases symmetrically in magnitude, as $\Delta$ increases in its range of dominance. Thus it is interesting to note that, at $q = \Delta^{-1}$, the proposed class of estimators is unbiased with largest efficiency and hence in the vicinity of $q = \Delta^{-1}$ also, the proposed class not only renders the massive gain in efficiency but also it is marginally biased in comparison of MMSE estimator. This implies that $q$ plays an important role in the proposed class of estimators. The following figure illustrates the discussion.

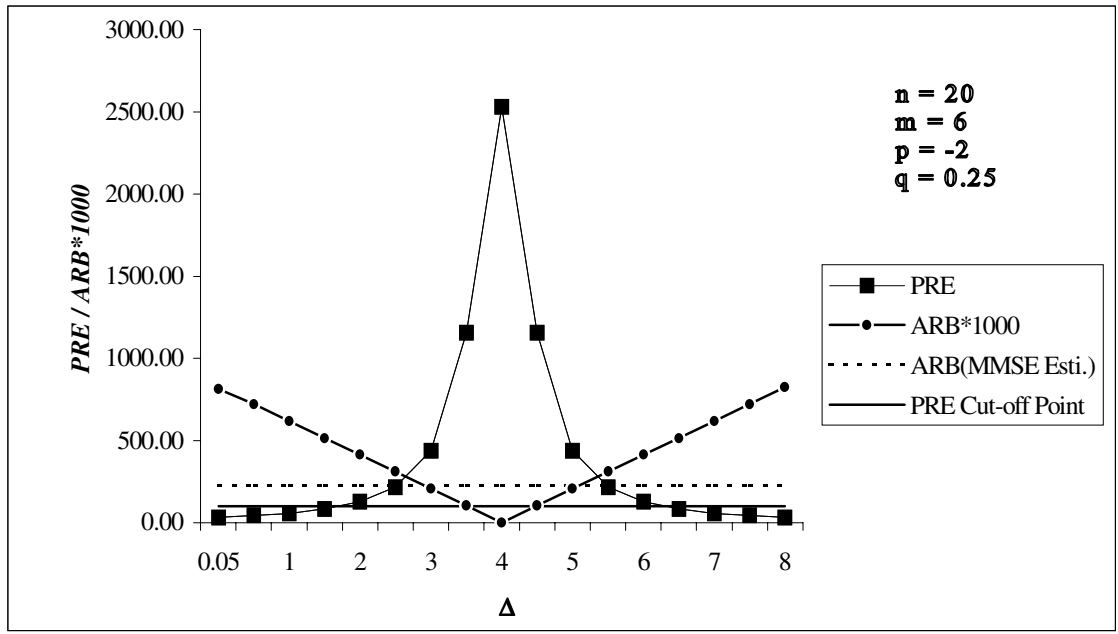

**Figure 3.1**

The effect of change in censored sample size $m$ is also a matter of great interest. For fixed $p$, $q$ and $\Delta$, the gain in relative efficiency diminishes, and ARB also decreases, with increment in $m$. Moreover, it appears that to get better estimators in the class, the value of $w(p)$ should be as small as possible in the interval (0,1]. Thus, to choose $p$ one should not consider the smaller values of $w(p)$ in isolation, but also the wider length of the interval of $\Delta$.

## 4. MODIFIED CLASS OF SHRINKAGE ESTIMATORS AND ITS PROPERTIES

The proposed class of estimators $\hat{\beta}_{(p,q)}$ is not uniformly better than $\hat{\beta}$. It will be better if $\beta_1$ and $\beta_2$ are in the vicinity of true value $\beta$. Thus, the centre of the guessed interval $(\beta_1 + \beta_2)/2$ is of much importance in this case. If we partially violate this, i.e., only the centre of the guessed interval is not of much importance, but the end points of the interval $\beta_1$ and $\beta_2$ are itself equally important then we can propose a new class of shrinkage estimators for the shape parameter $\beta$ by using the suggested class $\hat{\beta}_{(p,q)}$ as

$$\tilde{\beta}_{(p,q)} = \begin{cases} \beta_1 & , \text{ if } t > [(h-2)/\beta_1] \\ \left(\dfrac{h-2}{t}\right) w(p) + q\left(\dfrac{\beta_1+\beta_2}{2}\right)\{1-w(p)\} & , \text{ if } [(h-2)/\beta_2] \leq t \leq [(h-2)/\beta_1] \\ \beta_2 & , \text{ if } t < [(h-2)/\beta_2] \end{cases} \quad (4.1)$$

which has

$$\text{Bias}\{\tilde{\beta}_{(p,q)}\} = \beta\left[\Delta_1\left\{1 - I\left(\eta_1, \frac{h}{2}\right)\right\} + w(p)\left\{I\left(\eta_1, \frac{h}{2}-1\right) - I\left(\eta_2, \frac{h}{2}-1\right)\right\} \right. \\ \left. + q\Delta\{1-w(p)\}\left\{I\left(\eta_1, \frac{h}{2}\right) - I\left(\eta_2, \frac{h}{2}\right)\right\} + \Delta_2 I\left(\eta_2, \frac{h}{2}\right) - 1\right] \quad (4.2)$$

and

$$\text{MSE}\{\tilde{\beta}_{(p,q)}\} = \beta^2\left[(\Delta_1 - 1)^2 - \Delta_1(\Delta_1 - 2) I\left(\eta_1, \frac{h}{2}\right) + \Delta_2(\Delta_2 - 2) I\left(\eta_2, \frac{h}{2}\right) \right. \\ + \{w(p)\}^2 \left(\frac{h-2}{h-4}\right)\left\{I\left(\eta_1, \frac{h}{2}-2\right) - I\left(\eta_2, \frac{h}{2}-2\right)\right\} \\ + q\Delta\{1-w(p)\}\left\{I\left(\eta_1, \frac{h}{2}\right) - I\left(\eta_2, \frac{h}{2}\right)\right\}\{q\Delta\{1-w(p)\} - 2\} \\ \left. + 2w(p)\left\{I\left(\eta_1, \frac{h}{2}-1\right) - I\left(\eta_2, \frac{h}{2}-1\right)\right\}\{q\Delta\{1-w(p)\}-1\}\right] \quad (4.3)$$

where $\eta_1 = \left(\dfrac{h}{2} - 1\right)\Delta_1^{-1}$, $\eta_2 = \left(\dfrac{h}{2} - 1\right)\Delta_2^{-1}$ and $I(\eta, \omega) = \dfrac{1}{\Gamma(\omega)}\int_0^\eta e^{-u} u^{\omega-1} du$.

This modified class of shrinkage estimators is proposed in accordance with Rao(1973) and it seems to be more realistic than the previous one as it deals with the case where the whole interval is taken as apriori information.

## 5. NUMERICAL ILLUSTRATIONS

The percent relative efficiency of the proposed estimator $\tilde{\beta}_{(p,q)}$ with respect to MMSE estimator $\hat{\beta}_m$ has been defined as

$$\text{PRE}\{\tilde{\beta}_{(p,q)}, \hat{\beta}_m\} = \frac{\text{MSE}\{\hat{\beta}_m\}}{\text{MSE}\{\tilde{\beta}_{(p,q)}\}} \times 100 \qquad (5.1)$$

and it is obtained for $n = 20$ and different values of $p$, $q$, $m$, $\Delta_1$ and $\Delta_2$ (or $\Delta$). The findings are summarised in Table 5.1 with corresponding values of $h$ and $w(p)$.

**Table 5.1**

PREs of proposed estimator $\tilde{\beta}_{(p,q)}$ with respect to MMSE estimator $\hat{\beta}_m$

| | | | | | | | | | | | |
|---|---|---|---|---|---|---|---|---|---|---|---|
| | | | \multicolumn{8}{c|}{$n = 20$} | | |
| | | $p \rightarrow$ | \multicolumn{4}{c|}{-1} | | | | \multicolumn{4}{c|}{1} |
| $q\downarrow$ | | | $m \rightarrow$ | 6 | 8 | 10 | 12 | 6 | 8 | 10 | 12 |
| | $\Delta_1\downarrow$ | $\Delta_2\downarrow$ | $h \rightarrow$ | 10.8519 | 15.6740 | 20.8442 | 26.4026 | 10.8519 | 15.6740 | 20.8442 | 26.4026 |
| | | | $\Delta\downarrow\ w(p)\rightarrow$ | 0.7739 | 0.8537 | 0.8939 | 0.9180 | 0.6888 | 0.7737 | 0.8251 | 0.8779 |
| 0.25 | 0.2 | 0.3 | 0.25 | 50.80 | 41.39 | 34.91 | 30.59 | 49.84 | 40.10 | 34.66 | 31.15 |
| | 0.4 | 0.6 | 0.50 | 117.60 | 81.01 | 67.45 | 63.17 | 113.90 | 79.57 | 65.63 | 61.55 |
| | 0.6 | 0.9 | 0.75 | 261.72 | 227.42 | 203.08 | 172.06 | 227.59 | 191.97 | 172.31 | 156.69 |
| | 0.8 | 1.2 | 1.00 | 548.60 | 426.98 | 342.54 | 286.06 | 454.93 | 355.31 | 293.42 | 262.79 |
| | 1.0 | 1.5 | 1.25 | 649.95 | 470.44 | 375.91 | 314.98 | 636.21 | 504.49 | 427.74 | 353.74 |
| | 1.2 | 1.8 | 1.50 | 268.31 | 189.82 | 150.17 | 125.21 | 286.06 | 210.91 | 168.38 | 135.01 |
| | 1.5 | 2.0 | 1.75 | 80.46 | 53.66 | 39.90 | 31.38 | 82.35 | 55.10 | 40.79 | 31.74 |
| 0.50 | 0.2 | 0.3 | 0.25 | 50.84 | 41.32 | 34.76 | 30.39 | 49.90 | 40.03 | 34.45 | 30.87 |
| | 0.4 | 0.6 | 0.50 | 120.81 | 82.01 | 67.97 | 63.49 | 118.31 | 81.13 | 66.48 | 62.03 |
| | 0.6 | 0.9 | 0.75 | 298.17 | 253.12 | 221.74 | 184.38 | 271.73 | 225.47 | 198.40 | 173.57 |
| | 0.8 | 1.2 | 1.00 | 642.86 | 473.19 | 368.65 | 303.15 | 583.65 | 433.16 | 344.05 | 292.64 |
| | 1.0 | 1.5 | 1.25 | 626.09 | 435.87 | 345.16 | 289.53 | 658.77 | 481.87 | 390.95 | 317.87 |
| | 1.2 | 1.8 | 1.50 | 247.90 | 175.97 | 140.57 | 118.43 | 264.16 | 191.09 | 152.66 | 124.73 |
| | 1.5 | 2.0 | 1.75 | 78.41 | 52.66 | 39.39 | 31.11 | 79.96 | 53.72 | 40.02 | 31.36 |
| 0.75 | 0.2 | 0.3 | 0.25 | 50.89 | 41.24 | 34.60 | 30.19 | 49.97 | 39.95 | 34.23 | 30.59 |
| | 0.4 | 0.6 | 0.50 | 124.02 | 83.01 | 68.50 | 63.81 | 122.74 | 82.68 | 67.32 | 62.50 |
| | 0.6 | 0.9 | 0.75 | 339.92 | 282.24 | 242.46 | 197.73 | 325.66 | 266.36 | 229.58 | 192.68 |
| | 0.8 | 1.2 | 1.00 | 723.50 | 510.42 | 389.34 | 316.87 | 710.96 | 504.67 | 388.35 | 317.53 |
| | 1.0 | 1.5 | 1.25 | 566.19 | 392.47 | 312.16 | 263.77 | 597.64 | 421.61 | 337.17 | 278.26 |
| | 1.2 | 1.8 | 1.50 | 224.67 | 161.95 | 131.14 | 111.81 | 233.41 | 169.19 | 136.65 | 114.63 |
| | 1.5 | 2.0 | 1.75 | 76.05 | 51.59 | 38.85 | 30.83 | 76.93 | 52.14 | 39.17 | 30.95 |

*Table 5.1 continued ...*

| q↓ | Δ₁↓ | Δ₂↓ | p → m → h → Δ↓ w(p)→ | -2 6 10.8519 0.7739 | 8 15.6740 0.8537 | 10 20.8442 0.8939 | 12 26.4026 0.9180 | 2 6 10.8519 0.6888 | 8 15.6740 0.7737 | 10 20.8442 0.8251 | 12 26.4026 0.8779 |
|---|---|---|---|---|---|---|---|---|---|---|---|
| 0.25 | 0.2 | 0.3 | 0.25 | 46.04 | 34.18 | 30.92 | 30.53 | 46.77 | 34.81 | 30.96 | 31.23 |
|  | 0.4 | 0.6 | 0.50 | 92.48 | 72.59 | 59.44 | 53.42 | 98.00 | 73.36 | 59.48 | 54.88 |
|  | 0.6 | 0.9 | 0.75 | 106.83 | 95.44 | 92.75 | 90.11 | 128.68 | 102.24 | 93.16 | 100.45 |
|  | 0.8 | 1.2 | 1.00 | 145.02 | 131.16 | 126.15 | 122.15 | 191.47 | 145.23 | 126.97 | 144.22 |
|  | 1.0 | 1.5 | 1.25 | 220.29 | 243.10 | 282.54 | 320.74 | 305.32 | 273.81 | 284.60 | 368.42 |
|  | 1.2 | 1.8 | 1.50 | 208.14 | 211.32 | 202.36 | 179.81 | 250.20 | 220.57 | 202.56 | 175.49 |
|  | 1.5 | 2.0 | 1.75 | 82.08 | 57.89 | 43.07 | 33.36 | 84.21 | 57.95 | 43.06 | 33.12 |
| 0.50 | 0.2 | 0.3 | 0.25 | 46.28 | 34.31 | 30.86 | 30.24 | 46.95 | 34.91 | 30.90 | 30.87 |
|  | 0.4 | 0.6 | 0.50 | 103.18 | 76.82 | 61.54 | 54.80 | 107.21 | 77.31 | 61.57 | 56.08 |
|  | 0.6 | 0.9 | 0.75 | 157.81 | 135.64 | 127.02 | 118.59 | 181.60 | 142.94 | 127.44 | 128.23 |
|  | 0.8 | 1.2 | 1.00 | 267.16 | 228.67 | 207.62 | 190.69 | 331.58 | 246.71 | 208.58 | 212.20 |
|  | 1.0 | 1.5 | 1.25 | 445.44 | 443.06 | 448.55 | 438.38 | 541.60 | 467.49 | 449.42 | 432.21 |
|  | 1.2 | 1.8 | 1.50 | 289.70 | 240.03 | 198.56 | 163.98 | 298.93 | 238.16 | 198.30 | 156.40 |
|  | 1.5 | 2.0 | 1.75 | 84.92 | 57.28 | 42.13 | 32.67 | 84.44 | 57.03 | 42.12 | 32.44 |
| 0.75 | 0.2 | 0.3 | 0.25 | 46.50 | 34.43 | 30.78 | 29.92 | 47.13 | 34.99 | 30.82 | 30.50 |
|  | 0.4 | 0.6 | 0.50 | 114.64 | 81.04 | 63.59 | 56.13 | 116.87 | 81.23 | 63.61 | 57.24 |
|  | 0.6 | 0.9 | 0.75 | 247.11 | 202.90 | 181.31 | 160.85 | 266.60 | 209.00 | 181.65 | 167.34 |
|  | 0.8 | 1.2 | 1.00 | 543.26 | 418.40 | 345.15 | 293.90 | 596.79 | 430.93 | 345.67 | 302.22 |
|  | 1.0 | 1.5 | 1.25 | 704.42 | 541.77 | 447.06 | 381.03 | 696.36 | 532.12 | 446.25 | 358.48 |
|  | 1.2 | 1.8 | 1.50 | 280.39 | 203.46 | 160.74 | 132.95 | 269.47 | 199.82 | 160.55 | 129.07 |
|  | 1.5 | 2.0 | 1.75 | 81.39 | 54.49 | 40.40 | 31.66 | 80.35 | 54.26 | 40.39 | 31.52 |

It has been observed from Table 5.1 that likewise $\hat{\beta}_{(p,q)}$ the PRE of $\tilde{\beta}_{(p,q)}$ with respect to $\hat{\beta}_m$ decreases as censoring fraction ($m/n$) increases. For fixed $m$, $p$ and $q$ the relative efficiency increases up to a certain point of $\Delta$, procures its maximum at this point and then starts decreasing as $\Delta$ increases. It seems from the expression in (4.3) that the point of maximum efficiency may be a point where either any one of the following holds or any two of the following holds or all the following three holds-

(i) the lower end point of the guessed interval, i.e., $\beta_1$ coincides exactly with the true value $\beta$, i.e., $\Delta_1 = 1$.

(ii) the upper end point of the guessed interval, i.e., $\beta_2$ departs exactly two times from the true value $\beta$, i.e., $\Delta_2 = 2$.

(iii) $\Delta = q^{-1}$

This leads to say that on contrary to $\hat{\beta}_{(p,q)}$, there is much importance of $\Delta_1$ and $\Delta_2$ in addition to $\Delta$. The discussion is also supported by the illustrations in Table 5.1. As well, the range of dominance of average departure $\Delta$ is smaller than that is obtained for $\hat{\beta}_{(p,q)}$ but this does not humiliate the merit of $\tilde{\beta}_{(p,q)}$ because still the range of dominance of $\Delta$ is enough wider.

## 6. CONCLUSION AND RECOMMENDATIONS

It has been seen that the suggested classes of shrunken estimators have considerable gain in efficiency for a number of choices of scalars comprehend in it, particularly for heavily censored samples, i.e., for small *m*. Even for buoyantly censored samples, i.e., for large *m*, so far as the proper selection of scalars is concerned, some of the estimators from the suggested classes of shrinkage estimators are more efficient than the MMSE estimators subject to certain conditions. Accordingly, even if the experimenter has less confidence in the guessed interval $(\beta_1, \beta_2)$ of *β*, the efficiency of the suggested classes of shrinkage estimators can be increased considerably by choosing the scalars *p* and *q* appropriately.

While dealing with the suggested class of shrunken estimators $\hat{\beta}_{(p,q)}$ it is recommended that one should not consider the substantial gain in efficiency in isolation, but also the wider range of dominance of $\Delta$, because enough flexible range of dominance of $\Delta$ will leads to increase the possibility of getting better estimators from the proposed class. Thus it is recommended to use the proposed class of shrunken estimators in practice.